\documentclass[11pt]{amsart}
\pdfoutput=1

\usepackage{amssymb}
\usepackage{microtype}
\usepackage[numbers]{natbib}
\usepackage[pdftitle={$L^p$ change of variables inequalities on
manifolds},pdfauthor={Ari Stern},pdfsubject={MSC 2010: Primary 58C35;
Secondary 46B28, 58C40},bookmarksopen=false]{hyperref}

\theoremstyle{plain} 
\newtheorem{theorem}{Theorem}
\newtheorem{corollary}[theorem]{Corollary}
\theoremstyle{definition}
\theoremstyle{remark}
\newtheorem*{remark}{Remark}
\newtheorem*{acknowledgments}{Acknowledgments}

\begin{document}
\title[$ L ^p $ change of variables inequalities on
manifolds]{$\boldsymbol{L ^p}$ change of variables inequalities\\
  on manifolds}

\author{Ari Stern}

\address{%
Department of Mathematics\\
University of California, San Diego\\
9500 Gilman Dr \#0112\\
La Jolla CA 92093-0112}

\email{astern@math.ucsd.edu}

\begin{abstract}
  We prove two-sided inequalities for the $ L ^p $-norm of a
  pushforward or pullback (with respect to an orientation-preserving
  diffeomorphism) on oriented volume and Riemannian manifolds.  For a
  function or density on a volume manifold, these bounds depend only
  on the Jacobian determinant, which arises through the change of
  variables theorem.  For an arbitrary differential form on a
  Riemannian manifold, however, these bounds are shown to depend on
  more general ``spectral'' properties of the diffeomorphism, using an
  appropriately defined notion of singular values.  These spectral
  terms generalize the Jacobian determinant, which is recovered in the
  special cases of functions and densities (i.e., bottom and top
  forms).
\end{abstract}

\subjclass[2010]{Primary 58C35; Secondary 46B28, 58C40}

\keywords{Jacobian determinant, change of variables formula,
  differential forms, pushforward, pullback, singular values}

\date{November 28, 2011}

\maketitle

\section{Introduction}

\subsection{Motivation}
One of the most important tools in integral calculus (both classically
in $\mathbb{R}^n$ and on manifolds) is the change of variables
theorem, which states that the integral of a density is invariant
under pushforward and pullback by orientation-preserving
diffeomorphisms.  A closely related object is the Jacobian
determinant, which describes how the volume form changes with respect
to pushforwards and pullbacks.

This paper is motivated by the following, natural question:
\begin{quote}
  How do pushforwards and pullbacks affect the $ L ^p $-norm of a
  function or density, on an oriented volume manifold; or that of a
  differential form, on an oriented Riemannian manifold?
\end{quote}
We prove two-sided inequalities for each of these cases, showing that
the norms of the pushforward and pullback are controlled by
``spectral'' properties of the diffeomorphism, with respect to an
appropriately defined notion of singular values.  In the case of a
function or density, the bounds simply depend on the Jacobian
determinant of the diffeomorphism and that of its inverse; these can
be thought of as, respectively, the product of the singular values and
that of their reciprocals.  For a general differential form, though,
we encounter Jacobian determinant-like products that combine {\em
  both} singular values and reciprocal singular values.  In the $ L ^p
$-norm, these products are also shown to involve the conjugate
exponent $q$, which satisfies $ 1/p + 1/q = 1 $.

\subsection{Main results}
The main results of this paper are summarized in the following
theorems, which will be proved in the subsequent sections.

\begin{theorem}[smooth functions]
  \label{thm:scalarPushforward}
  Let $M$ and $N$ be oriented, $n$-dimensional manifolds with volume
  forms $ \mu _M $ and $ \mu _N $, respectively, and let $ \varphi
  \colon M \rightarrow N $ be an orientation-preserving diffeomorphism
  with Jacobian determinant $ J \left( \mu _M , \mu _N \right) \varphi
  \in C ^\infty (M) $.  Then, for any function $ u \in C ^\infty (M) $
  with compact support $ \operatorname{supp} u $,
  \begin{multline*}
    \bigl\lVert \mathbf{1} _{\operatorname{supp} u } \left[ J \left( \mu
        _M , \mu _N \right) \varphi \right] ^{-1/p} \bigr\rVert_\infty
    ^{-1} \left\lVert u \right\rVert _p \\
    \leq \left\lVert u \circ \varphi ^{-1} \right\rVert _p \leq
    \bigl\lVert \mathbf{1} _{ \operatorname{supp} u } \left[ J \left(
        \mu _M , \mu _N \right) \varphi \right] ^{ 1/p} \bigr\rVert
    _\infty \left\lVert u \right\rVert _p ,
  \end{multline*}
  for all $ p \in \left[ 1, \infty \right] $.
\end{theorem}

If the Jacobian determinant is bounded uniformly on all of $M$ (for
example, if $M$ and $N$ are compact), then this immediately yields a
uniform inequality,
\begin{equation*}
  \bigl\lVert \left[ J \left( \mu _M , \mu _N \right) \varphi \right]
  ^{-1/p} \bigr\rVert_\infty
  ^{-1} \left\lVert u \right\rVert _p
  \leq \left\lVert u \circ \varphi ^{-1} \right\rVert _p \leq
  \bigl\lVert \left[ J \left( \mu _M , \mu _N \right) \varphi \right]
  ^{ 1/p} \bigr\rVert _\infty \left\lVert u \right\rVert _p ,
\end{equation*}
which holds for all compactly supported $ u \in C^\infty (M) $.  This
implies that the map $ u \mapsto u \circ \varphi ^{-1} $ is bounded,
and because smooth functions with compact support form a dense subset
of $ L ^p (M) $, we can therefore extend this to the whole space.
Hence, \autoref{thm:scalarPushforward} has the following corollary.

\begin{corollary}[$ L ^p $ functions]
  Let $M$ and $N$ be oriented, $n$-dimensional manifolds with volume
  forms $ \mu _M $ and $ \mu _N $, respectively, and let $ \varphi
  \colon M \rightarrow N $ be an orientation-preserving diffeomorphism
  with Jacobian determinant $ J \left( \mu _M , \mu _N \right) \varphi
  \in C ^\infty (M) $.  If the Jacobian determinant is bounded
  uniformly on $M$, then for any $ u \in L ^p (M) $,
  \begin{equation*}
    \bigl\lVert \left[ J \left( \mu _M , \mu _N \right) \varphi \right]
    ^{-1/p} \bigr\rVert_\infty
    ^{-1} \left\lVert u \right\rVert _p
    \leq \left\lVert u \circ \varphi ^{-1} \right\rVert _p \leq
    \bigl\lVert \left[ J \left( \mu _M , \mu _N \right) \varphi \right]
    ^{ 1/p} \bigr\rVert _\infty \left\lVert u \right\rVert _p ,
  \end{equation*}
  for all $ p \in \left[ 1, \infty \right] $.
\end{corollary}

There is an analogous result for densities (i.e., $n$-forms) on a
volume manifold. (Since densities are the Hodge dual of functions, and
the $ L ^p $ and $ L ^q $ function spaces are dual to one another, it
is perhaps not too surprising that the conjugate exponent $q$ plays an
important role here.)

\begin{theorem}[smooth densities]
  \label{thm:densityPushforward}
  Let $M$ and $N$ be oriented, $n$-dimensional manifolds with volume
  forms $ \mu _M $ and $ \mu _N $, respectively, and let $ \varphi
  \colon M \rightarrow N $ be an orientation-preserving diffeomorphism
  with Jacobian determinant $ J \left( \mu _M , \mu _N \right) \varphi
  \in C ^\infty (M) $. Then, for any smooth density $ u \mu _M \in
  \Omega ^n (M) $, where $ u \in C^\infty (M) $ has compact support $
  \operatorname{supp} u = \operatorname{supp} u \mu _M $,
  \begin{multline*}
    \bigl\lVert \mathbf{1} _{ \operatorname{supp} u \mu _M }\left[ J
      \left( \mu _M , \mu _N \right) \varphi \right] ^{ 1/q }
    \bigr\rVert _\infty ^{-1}
    \left\lVert u \mu _M \right\rVert _p  \\
    \leq \left\lVert \varphi _\ast \left( u \mu _M \right) \right\rVert
    _p \leq \bigl\lVert \mathbf{1} _{ \operatorname{supp} u \mu _M }
    \left[ J \left( \mu _M , \mu _N \right) \varphi \right] ^{ -1/q }
    \bigr\rVert _\infty \left\lVert u \mu _M \right\rVert _p ,
  \end{multline*}
  for all $ p , q \in \left[ 1, \infty \right] $ such that $ 1/p + 1/q
  = 1 $.
\end{theorem}

\begin{corollary}[$ L ^p $ densities]
  Let $M$ and $N$ be oriented, $n$-dimensional manifolds with volume
  forms $ \mu _M $ and $ \mu _N $, respectively, and let $ \varphi
  \colon M \rightarrow N $ be an orientation-preserving diffeomorphism
  with Jacobian determinant $ J \left( \mu _M , \mu _N \right) \varphi
  \in C ^\infty (M) $.  If the reciprocal $ \left[ J \left( \mu _M ,
      \mu _N \right) \varphi \right] ^{-1} \mkern-3mu $ is bounded
  uniformly on $M$, then for any density $ u \mu _M \in L ^p \Omega ^n
  (M) $, where $ u \in L ^p (M) $,
  \begin{multline*}
    \bigl\lVert \left[ J \left( \mu _M , \mu _N \right) \varphi
    \right] ^{ 1/q } \bigr\rVert _\infty ^{-1}
    \left\lVert u \mu _M \right\rVert _p  \\
    \leq \left\lVert \varphi _\ast \left( u \mu _M \right)
    \right\rVert _p \leq \bigl\lVert \left[ J \left( \mu _M , \mu _N
      \right) \varphi \right] ^{ -1/q } \bigr\rVert _\infty
    \left\lVert u \mu _M \right\rVert _p ,
  \end{multline*}
  for all $ p , q \in \left[ 1, \infty \right] $ such that $ 1/p + 1/q
  = 1 $.
\end{corollary}

When $M$ and $N$ are Riemannian manifolds, however, we show that it is
possible to obtain a much more general family of inequalities, which
hold for arbitrary $k$-forms on $M$, where $ k = 0 , 1, \ldots, n $.
In the special cases $ k = 0 $ and $ k = n $, these inequalities are
shown to recover the previous results for functions and densities,
respectively.  In order to state and prove these more general results,
we will introduce a novel extension of the \emph{singular values} of a
mapping.  Whereas singular values are traditionally defined only for
linear maps between Euclidean vector spaces, we show that they can
also be defined intrinsically for diffeomorphisms between Riemannian
manifolds.

\begin{theorem}[smooth $k$-forms]
  \label{thm:kformPushforward}
  Let $ \left( M, g _M \right) $ and $ \left( N, g _N \right) $ be
  oriented, $n$-dimensional Riemannian manifolds, and let $ \varphi
  \colon M \rightarrow N $ be an orientation-preserving diffeomorphism
  with singular values $ \alpha _1 (x) \geq \cdots \geq \alpha _n (x)
  > 0 $ at each $ x \in M $.  Then, for any smooth $k$-form $ \omega
  \in \Omega ^k (M) $, $ k = 0, \ldots, n $, with compact support $
  \operatorname{supp} \omega $,
  \begin{multline*}
    \bigl\lVert \mathbf{1} _{\operatorname{supp} \omega } \left( \alpha
      _1 \cdots \alpha _k \right) ^{ 1/q } \left( \alpha _{ k + 1 }
      \cdots \alpha _n \right) ^{-1/p} \bigr\rVert _\infty ^{-1}
    \left\lVert \omega
    \right\rVert _p  \\
    \leq \left\lVert \varphi _\ast \omega \right\rVert _p \leq
    \bigl\lVert \mathbf{1} _{ \operatorname{supp} \omega } \left( \alpha
      _1 \cdots \alpha _{ n - k } \right) ^{1/p} \left( \alpha _{ n - k
        + 1 } \cdots \alpha _n \right) ^{ -1/q }\bigr\rVert _\infty
    \left\lVert \omega \right\rVert _p ,
  \end{multline*}
  for all $ p , q \in \left[ 1, \infty \right] $ such that $ 1/p + 1/q
  = 1 $.
\end{theorem}

\begin{corollary}[$ L ^p $ $k$-forms]
  Let $ \left( M, g _M \right) $ and $ \left( N, g _N \right) $ be
  oriented, $n$-dimen\-sional Riemannian manifolds, and let $ \varphi
  \colon M \rightarrow N $ be an orientation-preserving diffeomorphism
  with singular values $ \alpha _1 (x) \geq \cdots \geq \alpha _n (x)
  > 0 $ at each $ x \in M $.  Given $ p , q \in \left[ 1, \infty
  \right] $ such that $ 1/p + 1/q = 1 $, and some $ k = 0 , \ldots, n
  $, suppose that the product $ \left( \alpha _1 \cdots \alpha _{ n -
      k } \right) ^{1/p} \left( \alpha _{ n - k + 1 } \cdots \alpha _n
  \right) ^{ -1/q } $ is bounded uniformly on $M$.  Then, for any $
  \omega \in L ^p \Omega ^k (M) $,
  \begin{multline*}
    \bigl\lVert \left( \alpha _1 \cdots \alpha _k \right) ^{ 1/q }
    \left( \alpha _{ k + 1 } \cdots \alpha _n \right) ^{-1/p}
    \bigr\rVert _\infty ^{-1} \left\lVert \omega
    \right\rVert _p  \\
    \leq \left\lVert \varphi _\ast \omega \right\rVert _p \leq
    \bigl\lVert \left( \alpha _1 \cdots \alpha _{ n - k } \right)
    ^{1/p} \left( \alpha _{ n - k + 1 } \cdots \alpha _n \right) ^{
      -1/q }\bigr\rVert _\infty \left\lVert \omega \right\rVert _p .
  \end{multline*}
\end{corollary}

\subsection{Organization of the paper}
We begin, in the next section, by proving
\autoref{thm:scalarPushforward} and \autoref{thm:densityPushforward},
which apply to functions and densities, respectively, on oriented
volume manifolds.

In the subsequent section, we turn to the case of differential
$k$-forms on a Riemannian manifold.  Unlike the previous inequalities
for volume manifolds, whose main ingredients are the change of
variables theorem and H\"older's inequality, the Riemannian case
requires completely new analytical tools, which we introduce and
develop along the way.  In particular, the proof of
\autoref{thm:kformPushforward} depends crucially on the generalized
definition of singular values for a diffeomorphism between Riemannian
manifolds.  (Note that this is distinct from the usual notion of
\emph{spectrum} for a Riemannian manifold, which typically refers to
eigenvalues of the Laplacian \citep{CrPuRa2001}.) This theorem, along
with the new techniques required to state and prove it, are the most
significant contributions of this paper.

The proof of \autoref{thm:kformPushforward} also depends on some facts
from multilinear algebra, relating the singular values of a linear
operator to the spectral norm of its induced map on alternating
tensors. We provide supplementary technical details in
\autoref{app:multilinear}.

Finally, we mention that although these results only apply to
pushforwards, as stated above, it is trivial to apply them to
pullbacks as well.  Using the definition $ \varphi _\ast = \left(
  \varphi ^{-1} \right) ^\ast $, simply replace $ \varphi $ by $
\varphi ^{-1} $ above, along with its corresponding Jacobian
determinant and singular values.  For completeness (and for the
convenience of the reader), the pullback versions of these
inequalities are stated in \autoref{app:pullback}.

\section{Change of variables on volume manifolds}

\subsection{Smooth functions with compact support} Let $M$ and $N$ be
oriented, $n$-dimensional manifolds, with volume forms $ \mu
_M $ and $ \mu _N $, respectively.  If $ \varphi \colon M \rightarrow
N $ is an orientation-preserving diffeomorphism, recall that there
exists a function $ J \left( \mu _M , \mu _N \right) \varphi \in C
^\infty (M) $, called the {\em Jacobian determinant}, such that $
\varphi ^\ast \mu _N = \left[ J \left( \mu _M , \mu _N \right) \varphi
\right] \mu _M $ (see, e.g., \citet*{AbMaRa1988}).  For any function $
v \in C ^\infty (N) $ with compact support, this implies the familiar
change of variables formula
\begin{equation*}
  \int _N v \mu _N = \int _M \varphi ^\ast \left( v \mu _N \right) =
  \int _M \left( v \circ \varphi \right) \left[ J \left( \mu _M , \mu _N
    \right) \varphi \right] \mu _M .
\end{equation*}
Now, let $ u \in C^\infty (M) $ be a function with compact support
(denoted by $ \operatorname{supp} u $).  Then, for any $ p \in
\left[1, \infty \right) $, it follows that
\begin{align*}
  \int _N \left\lvert u \circ \varphi ^{-1} \right\rvert ^p \mu _N &=
  \int _M \left\lvert u \right\rvert ^p \left[ J \left( \mu _M, \mu _N
    \right) \varphi \right] \mu _M \\
  &\leq \left\lVert \mathbf{1} _{ \operatorname{supp} u } J \left( \mu
      _M , \mu _N \right) \varphi \right\rVert _\infty \int _M
  \left\lvert u \right\rvert ^p \mu _M ,
\end{align*}
by change of variables and H\"older's inequality.  This immediately
gives the $ L ^p $-norm upper bound
\begin{equation*}
  \left\lVert u \circ \varphi ^{-1} \right\rVert _p \leq \bigl\lVert
  \mathbf{1} _{ \operatorname{supp} u } \left[ J \left( \mu
      _M , \mu _N \right) \varphi \right] ^{ 1/p} \bigr\rVert
  _\infty \left\lVert u \right\rVert _p .
\end{equation*}
(Notice that this also holds for $ p = \infty $.)
To obtain the lower bound, we write
\begin{align*}
  \int _M \left\lvert u \right\rvert ^p \mu _M &= \int _M \left(
    \left\lvert u \circ \varphi ^{-1} \right\rvert ^p \circ \varphi
  \right) \mu _M \\
  &= \int _M \varphi ^\ast \left( \left\lvert u \circ \varphi
    \right\rvert ^{-1} \varphi _\ast \mu _M \right) \\
  &= \int _M \varphi ^\ast \left( \left\lvert u \circ \varphi
    \right\rvert ^{-1} \left[ J \left( \mu _N , \mu _M \right) \left(
        \varphi ^{-1} \right) \right] \mu _N \right)
\end{align*}
Since the Jacobian determinant satisfies the inverse identity
\begin{equation}
  \label{eqn:jacobianInverse}
  J \left( \mu _N , \mu _M \right) \left( \varphi ^{-1} \right) =   \left[J
    \left( \mu _M , \mu _N \right) \varphi \right] ^{-1} \circ
  \varphi ^{-1} ,
\end{equation}
it follows that
\begin{align*}
  \int _M \left\lvert u \right\rvert ^p \mu _M &= \int _M \left[ J
    \left( \mu _M , \mu _N \right) \varphi \right] ^{-1} \varphi ^\ast
  \left( \left\lvert u \circ \varphi ^{-1} \right\rvert ^p \mu _N
  \right) \\
  &\leq \bigl\lVert \mathbf{1} _{ \operatorname{supp} u } \left[ J
    \left( \mu _M , \mu _N \right) \varphi \right] ^{-1} \bigr\rVert
  _\infty \int _N \left\lvert u \circ \varphi ^{-1} \right\rvert ^p
  \mu _N,
\end{align*}
again using H\"older's inequality and change of variables.  Thus,
\begin{equation*}
  \left\lVert u \right\rVert _p \leq \bigl\lVert \mathbf{1} _{
    \operatorname{supp} u } \left[ J  \left( \mu _M , \mu _N \right)
    \varphi \right] ^{-1/p} \bigr\rVert_\infty \left\lVert u \circ
    \varphi ^{-1} \right\rVert _p ,
\end{equation*}
which rearranges to give the lower bound
\begin{equation*}
  \left\lVert u \circ \varphi ^{-1} \right\rVert _p  \geq \bigl\lVert
  \mathbf{1} _{\operatorname{supp} u }  \left[ J \left( \mu _M , \mu _N
    \right) \varphi \right] ^{-1/p} \bigr\rVert_\infty ^{-1}
  \left\lVert u \right\rVert _p .
\end{equation*}
In summary, we have now established the two-sided inequality,
\begin{multline*}
  \bigl\lVert \mathbf{1} _{\operatorname{supp} u } \left[ J \left( \mu
      _M , \mu _N \right) \varphi \right] ^{-1/p} \bigr\rVert_\infty
  ^{-1} \left\lVert u \right\rVert _p \\
  \leq \left\lVert u \circ \varphi ^{-1} \right\rVert _p \leq
  \bigl\lVert \mathbf{1} _{ \operatorname{supp} u } \left[ J \left(
      \mu _M , \mu _N \right) \varphi \right] ^{ 1/p} \bigr\rVert
  _\infty \left\lVert u \right\rVert _p ,
\end{multline*}
which completes the proof of \autoref{thm:scalarPushforward}. \qed

\begin{remark}
  If $\varphi$ is volume-preserving (at least on the support of $u$),
  then it also preserves the $ L ^p $-norm for all $p$.  Indeed, if $
  J \left( \mu _M , \mu _N \right) \varphi = 1 $, then the inequality
  simply becomes $ \left\lVert u \right\rVert _p \leq \left\lVert u
    \circ \varphi ^{-1} \right\rVert _p \leq \left\lVert u
  \right\rVert _p $, and thus $ \left\lVert u \circ \varphi ^{-1}
  \right\rVert _p = \left\lVert u \right\rVert _p $, as expected.
  This equality is also seen to hold for arbitrary diffeomorphisms
  $\varphi$ (i.e., not necessarily volume-preserving) when $ p =
  \infty $.
\end{remark}

\subsection{Smooth densities with compact support}
Any smooth density on $N$ with compact support can be written as $ v
\mu _N $, where $ v \in C ^\infty (N) $ is a compactly-supported
function.  Now, since the Hodge star operator $ \star _N $ is an
isometry, with $ \star _N \mu _N = 1 $, the pointwise norm of this
density satisfies $ \left\lvert v \mu _N \right\rvert = \left\lvert
  \star _N \left( v \mu _N \right) \right\rvert = \left\lvert v
\right\rvert $ (and likewise for $ \star _M $).  In particular, this
implies the identity
\begin{equation*}
  \left\lvert v \mu _N \right\rvert ^p \circ \varphi = \left\lvert v
  \right\rvert ^p \circ \varphi = \left\lvert v \circ \varphi
  \right\rvert ^p = \left\lvert \left( v \circ \varphi \right) \mu _M
  \right\rvert ^p ,
\end{equation*}
Therefore,
\begin{equation}
  \label{eqn:densityIdentity}
  \begin{split}
    \left\lvert \varphi ^\ast \left( v \mu _N \right) \right\rvert ^p
    &= \left\lvert \left( v \circ \varphi \right) \left[ J \left( \mu
          _M , \mu _N \right) \varphi \right] \mu _M \right\rvert ^p
    \\
    &= \left( \left\lvert v \mu _N \right\rvert ^p \circ \varphi
    \right) \left[ J \left( \mu _M , \mu _N \right) \varphi \right] ^p
    .
  \end{split}
\end{equation}

Now, consider the pushforward of the density $ u \mu _M $, where $ u
\in C^\infty (M) $ has compact support $ \operatorname{supp} u =
\operatorname{supp} u \mu _M $.  Using change of variables,
\begin{align*}
  \int _N \left\lvert \varphi _\ast \left( u \mu _M \right)
  \right\rvert ^p \mu _N &= \int _M \left( \left\lvert \varphi _\ast
      \left( u \mu _M \right) \right\rvert ^p \circ \varphi \right)
  \left[ J \left( \mu _M , \mu _N \right) \varphi \right] \mu _M \\
  &= \int _M \left\lvert u \mu _M \right\rvert ^p \left[ J \left( \mu
      _M , \mu _N \right) \varphi \right] ^{ 1 - p } \mu _M \\
  &\leq \bigl\lVert \mathbf{1} _{ \operatorname{supp} u \mu _M }
  \left[ J \left( \mu _M , \mu _N \right) \varphi \right] ^{ 1 - p }
  \bigr\rVert _\infty \int _M \left\lvert u \mu _M \right\rvert ^p \mu
  _M ,
\end{align*}
where the last two lines use the identity \eqref{eqn:densityIdentity}
and H\"older's inequality, respectively.  Taking the conjugate
exponent $q$ such that $ 1/p + 1/q = 1 $, observe that $ \left( 1 - p
\right) / p = - 1/q $.  Therefore, we can write the upper bound
\begin{equation*}
  \left\lVert  \varphi _\ast \left( u \mu _M \right) \right\rVert _p
  \leq \bigl\lVert \mathbf{1} _{ \operatorname{supp} u \mu _M }
  \left[ J \left( \mu _M , \mu _N \right) \varphi \right] ^{ -1/q }
  \bigr\rVert _\infty \left\lVert u \mu _M \right\rVert _p .
\end{equation*}
For the lower bound, we begin by using the identity
\eqref{eqn:densityIdentity} to write
\begin{equation*}
  \left\lvert u \mu _M \right\rvert ^p = \left\lvert \varphi ^\ast
    \varphi _\ast \left( u \mu _M \right) \right\rvert ^p = \left(
    \left\lvert \varphi _\ast \left( u \mu _M \right) \right\rvert ^p
    \circ \varphi \right) \left[ J \left( \mu _M , \mu _N \right)
    \varphi \right] ^p ,
\end{equation*} 
and thus
\begin{align*}
  \int _M \left\lvert u \mu _M \right\rvert ^p \mu _M &= \int _M
  \left( \left\lvert \varphi _\ast \left( u \mu _M \right)
    \right\rvert ^p \circ \varphi \right) \left[ J \left( \mu _M , \mu
      _N \right) \varphi \right] ^p \mu _M \\
  &= \int _M \left[ J \left( \mu _M , \mu _N \right) \varphi \right]
  ^p \varphi ^\ast \left( \left\lvert \varphi _\ast \left( u \mu _M
      \right) \right\rvert ^p \varphi _\ast \mu _M \right) \\
  &= \int _M \left[ J \left( \mu _M , \mu _N \right) \varphi \right]
  ^{p-1} \varphi ^\ast \left( \left\lvert \varphi _\ast \left( u \mu
        _M \right) \right\rvert ^p \mu _N \right) \\
  &\leq \bigl\lVert \mathbf{1} _{ \operatorname{supp} u \mu _M }
  \left[ J \left( \mu _M , \mu _N \right) \varphi \right] ^{p-1}
  \bigr\rVert _\infty \int _N \left\lvert \varphi _\ast \left( u \mu
      _M \right) \right\rvert ^p \mu _N ,
\end{align*}
by the Jacobian inverse identity \eqref{eqn:jacobianInverse},
H\"older's inequality, and change of variables.  Therefore,
\begin{equation*}
  \left\lVert u \mu _M \right\rVert _p \leq \bigl\lVert \mathbf{1} _{
    \operatorname{supp} u \mu _M }\left[ J \left( \mu _M , \mu _N
    \right) \varphi \right] ^{ 1/q } \bigr\rVert _\infty \left\lVert
    \varphi _\ast \left( u \mu _M \right) \right\rVert _p ,
\end{equation*}
so rearranging, we have the lower bound
\begin{equation*}
  \left\lVert \varphi _\ast \left( u \mu _M \right) \right\rVert _p
  \geq \bigl\lVert \mathbf{1} _{
    \operatorname{supp} u \mu _M }\left[ J \left( \mu _M , \mu _N
    \right) \varphi \right] ^{ 1/q } \bigr\rVert _\infty ^{-1}
  \left\lVert u \mu _M \right\rVert _p .
\end{equation*}
Hence, we have shown the two-sided inequality
\begin{multline*}
  \bigl\lVert \mathbf{1} _{ \operatorname{supp} u \mu _M }\left[ J
    \left( \mu _M , \mu _N \right) \varphi \right] ^{ 1/q }
  \bigr\rVert _\infty ^{-1}
  \left\lVert u \mu _M \right\rVert _p  \\
  \leq \left\lVert \varphi _\ast \left( u \mu _M \right) \right\rVert
  _p \leq \bigl\lVert \mathbf{1} _{ \operatorname{supp} u \mu _M }
  \left[ J \left( \mu _M , \mu _N \right) \varphi \right] ^{ -1/q }
  \bigr\rVert _\infty \left\lVert u \mu _M \right\rVert _p ,
\end{multline*}
which completes the proof of \autoref{thm:densityPushforward}. \qed

\begin{remark}
  For densities, the $ L ^p $-norm is preserved if either $\varphi$ is
  volume-preserving (at least on $ \operatorname{supp} u \mu _M $) or
  when $ p = 1 $, $ q = \infty $.  In these cases, the inequality
  becomes $ \left\lVert u \mu _M \right\rVert _p \leq \left\lVert
    \varphi _\ast \left( u \mu _M \right) \right\rVert _p \leq
  \left\lVert u \mu _M \right\rVert _p $, and thus $ \left\lVert
    \varphi _\ast \left( u \mu _M \right) \right\rVert _p =
  \left\lVert u \mu _M \right\rVert _p $.
\end{remark}

\section{Change of variables on Riemannian manifolds}

\subsection{Singular values of a diffeomorphism}
In order to generalize this result to differential forms, we suppose
now that $ \left( M , g _M \right) $ and $ \left( N, g _N \right) $
are oriented, $n$-dimensional Riemannian manifolds, with $ \mu _M $
and $ \mu _N $ denoting their respective Riemannian volume forms.  As
before, let $ \varphi \colon M \rightarrow N $ be an
orientation-preserving diffeomorphism.  Given a point $ x \in M $, let
$ \left\{ e _1 , \ldots , e _n \right\} $ be a positively-oriented, $
g _M $-orthonormal basis of the tangent space $ T _x M $, and let $
\left\{ f _1, \ldots, f _n \right\} $ be a positively-oriented, $ g _N
$-orthonormal basis of $ T _{ \varphi (x) } N $.  With respect to
these bases, the tangent map $ T _x \varphi \colon T _x M \rightarrow
T _{ \varphi (x) } N $ can be represented by an $ n \times n $ matrix
$ \Phi $.  Since $\varphi$ is a diffeomorphism, the matrix $\Phi $ has
$n$ positive singular values, which we write
\begin{equation*}
  \alpha _1 (x) \geq \cdots \geq \alpha _n (x) > 0 .
\end{equation*}
The singular values of $\Phi$ are orthogonally invariant, so they are
independent of the choice of orthonormal basis, and thus are an
intrinsic property of the diffeomorphism.  Therefore, we refer to
these as the {\em singular values of $\varphi$ at $x$}.  It follows
that the pullback of the volume form on $N$ is
\begin{equation*}
  \varphi ^\ast \mu _N = \left( \det \Phi \right) \mu _M = \left( \alpha _1
    \cdots \alpha _n \right) \mu _M ,
\end{equation*}
so the Jacobian determinant is simply the product of the singular
values $ J \left( \mu _M , \mu _N \right) \varphi = \alpha _1 \cdots
\alpha _n $.

Similarly, the inverse map $ T _{\varphi (x)} \left( \varphi ^{-1}
\right) \colon T _{ \varphi (x) } N \rightarrow T _x M $ is
represented by the inverse matrix $ \Phi ^{-1} $, whose singular
values are the reciprocals of those for $\Phi$.  Hence, we write the
singular values of $ \varphi ^{-1} $ at $ \varphi (x) $ as
\begin{equation*}
  \beta _1 \left( \varphi (x) \right) \geq \cdots \geq \beta _n \left(
    \varphi (x) \right) > 0 ,
\end{equation*}
which satisfy $ \beta _i \left( \varphi (x) \right) = \alpha _{ n - i
  + 1 } (x) ^{-1} $, i.e., $ \beta _i = \alpha _{ n - i + 1 } ^{-1}
\circ \varphi ^{-1} $, for $ i = 1, \ldots, n $.  Consequently, the
pushforward of the volume form on $M$ is 
\begin{equation*}
  \varphi _\ast \mu _M = \left( \det \Phi ^{-1} \right) \mu _N =
  \left( \beta _1 \cdots \beta _n \right) \mu _N .
\end{equation*}
Therefore, the Jacobian determinant is
\begin{multline*}
  J \left( \mu _N , \mu _M \right) \left( \varphi ^{-1} \right) \\
  = \beta _1 \cdots \beta _n = \left( \alpha _1 \cdots \alpha _n
  \right) ^{-1} \circ \varphi ^{-1} = \left[ J \left( \mu _M , \mu _N
    \right) \varphi \right] ^{-1} \circ \varphi ^{-1} ,
\end{multline*}
so we recover the usual identity \eqref{eqn:jacobianInverse}.

\begin{remark}
  Using the well-known ``minimax'' and ``maximin'' characterizations
  of singular values, it is also possible to write
  \begin{alignat*}{2}
    \alpha _i (x) &= \min _{ \substack{ S \subset \mathbb{R}^n  \\
        \dim S = n - i + 1 } } \max _{ 0 \neq X \in S } \frac{
      \left\lvert \Phi X \right\rvert }{ \left\lvert X \right\rvert }
    &&=  \max _{ \substack{ S \subset \mathbb{R}^n  \\
        \dim S = i } } \min _{ 0 \neq X \in S } \frac{ \left\lvert
        \Phi X \right\rvert }{ \left\lvert X \right\rvert } \\
    &= \min _{ \substack{ S \subset T _x M  \\
        \dim S = n - i + 1 } } \max _{ 0 \neq X \in S } \frac{
      \left\lvert T _x \varphi (X) \right\rvert }{ \left\lvert X
      \right\rvert } &&= \max
    _{ \substack{ S \subset T _x M   \\
        \dim S = i } } \min _{ 0 \neq X \in S } \frac{ \left\lvert T
        _x \varphi (X) \right\rvert }{ \left\lvert X \right\rvert }.
  \end{alignat*}
  This can be taken as an alternative, basis-independent definition
  for the singular values of a diffeomorphism, consistent with the
  previous one.  This also provides another way to see that $ \beta _i
  = \alpha _{ n - i + 1 } \circ \varphi ^{-1} $, since
  \begin{align*}
    \beta _i \left( \varphi (x) \right) &= \min _{ \substack{ S
        \subset T _{\varphi (x) } N \\ \dim S = n - i + 1 } } \max _{
      0 \neq Y \in S } \frac{ \left\lvert T _{ \varphi (x) } \left(
          \varphi ^{-1} \right) (Y) \right\rvert }{ \left\lvert Y
      \right\rvert } \\
    &= \min _{ \substack{ S \subset T _x M \\ \dim S = n - i + 1 } }
    \max _{ 0 \neq X \in S } \frac{ \left\lvert X \right\rvert }{
      \left\lvert T _x \varphi (X) \right\rvert } \\
    &= \biggl( \max _{ \substack{ S \subset T _x M \\ \dim S = n - i +
        1 } } \min _{ 0 \neq X \in S } \frac{ \left\lvert T _x \varphi
        (X) \right\rvert }{ \left\lvert X \right\rvert } \biggr) ^{-1}
    \\
    &= \alpha _{ n - i + 1 } (x) ^{-1} .
  \end{align*}
\end{remark}

\subsection{Pointwise inequalities for the spectral norm}
Now, given a smooth $k$-form $ \omega \in \Omega ^k (M) $ with compact
support, recall that the spectral norm of $\omega$ at a point $ x \in
M $ is defined by
\begin{equation*}
  \left\lvert \omega \right\rvert = \max _{ 0 \neq X _1, \ldots, X _k
    \in T _x M } \frac{ \left\lvert \omega \left( X _1, \ldots, X _k
      \right) \right\rvert }{ \left\lvert X _1 \right\rvert \cdots
    \left\lvert X _k \right\rvert }, 
\end{equation*}
where $ \left\lvert X _i \right\rvert = g _M \left( X _i , X _i
\right) ^{ 1/2} $ denotes the length of the tangent vector $ X _i $.
With this definition, the $ L ^p $-norm on $ \Omega ^k (M) $ is simply
$ \left\lVert \omega \right\rVert _p = \int _M \left\lvert \omega
\right\rvert ^p \mu _M $, as usual.  Likewise, for $ \eta \in \Omega
^k (N) $, at each point $ y \in N $ we have
\begin{equation*}
  \left\lvert \eta  \right\rvert = \max _{ 0 \neq Y _1, \ldots, Y _k
    \in T _y N } \frac{ \left\lvert \eta  \left( Y _1, \ldots, Y _k
      \right) \right\rvert }{ \left\lvert Y _1 \right\rvert \cdots
    \left\lvert Y _k \right\rvert }, 
\end{equation*}
where here $ \left\lvert Y _i \right\rvert = g _N \left( Y _i , Y _i
\right) ^{ 1/2 } $, and $ \left\lVert \eta \right\rVert _p = \int _N
\left\lvert \eta \right\rvert ^p \mu _N $.

To prove \autoref{thm:kformPushforward}, we begin by stating pointwise
bounds for the pullback and pushforward, in terms of the singular
values of $\varphi$ and $ \varphi ^{-1} $.  Since $\omega$ and $\eta$
are $k$-linear and totally antisymmetric at each point, it is
straightforward to show that
\begin{alignat}{2}
  \left\lvert \varphi ^\ast \eta \right\rvert &\leq \alpha _1 \cdots
  \alpha _k \left( \left\lvert \eta \right\rvert \circ \varphi \right)
  &&= \bigl[ \left( \beta _{ n - k + 1 } \cdots \beta _n \right) ^{-1}
    \left\lvert \eta \right\rvert \bigr] \circ
  \varphi \label{eqn:pointwisePullback}\\
  \left\lvert \varphi _\ast \omega \right\rvert &\leq \beta _1 \cdots
  \beta _k \left( \left\lvert \omega \right\rvert \circ \varphi ^{-1}
  \right) &&= \bigl[ \left( \alpha _{ n - k + 1 } \cdots \alpha _n
    \right) ^{-1} \left\lvert \omega \right\rvert \bigr] \circ
  \varphi ^{-1} \label{eqn:pointwisePushforward}.
\end{alignat}
That is, the pullback of a $k$-form is controlled by the product of
the $k$ largest singular values of $\varphi$, while the pushforward is
controlled by the product of the $k$ largest singular values of $
\varphi ^{-1} $.  Further discussion of these pointwise inequalities,
as well as some background and details on their derivation, is given
in \autoref{app:multilinear}.

\subsection{Change of variables for $k$-forms}
Using change of variables, the pointwise inequality
\eqref{eqn:pointwisePushforward} for the pushforward, and H\"older's
inequality,
\begin{align*}
  \int _N \left\lvert \varphi _\ast \omega \right\rvert ^p \mu _N &=
  \int _M \varphi ^\ast \left( \left\lvert \varphi _\ast \omega
    \right\rvert ^p \mu _N \right) \\
  &= \int _M \left( \left\lvert \varphi _\ast \omega
    \right\rvert ^p  \circ \varphi \right) \varphi ^\ast \mu _N \\
  &\leq \int _M \bigl[ \left( \alpha _{ n - k + 1 } \cdots \alpha _n
  \right) ^{-1} \left\lvert \omega \right\rvert \bigr] ^p \left(
    \alpha _1 \cdots \alpha _n \right) \mu _M \\
  &= \int _M \left( \alpha _1 \cdots \alpha _{ n - k } \right) \left(
    \alpha _{ n - k + 1 } \cdots \alpha _n \right) ^{ 1 - p }
  \left\lvert \omega \right\rvert ^p \mu _M \\
  &\leq \bigl\lVert \mathbf{1} _{ \operatorname{supp} \omega } \left(
    \alpha _1 \cdots \alpha _{ n - k } \right) \left( \alpha _{ n - k
      + 1 } \cdots \alpha _n \right) ^{ 1 - p } \bigr\rVert _\infty
  \int _M \left\lvert \omega \right\rvert ^p \mu _M .
\end{align*}
Hence, we immediately obtain the upper bound
\begin{equation*}
  \left\lVert \varphi _\ast \omega \right\rVert _p \leq \bigl\lVert
  \mathbf{1} _{ \operatorname{supp} \omega } \left( \alpha _1 \cdots
    \alpha _{ n - k } \right) ^{ 1/p } \left( \alpha _{ n - k + 1 } \cdots
    \alpha _n \right) ^{ -1/q } \bigr\rVert _\infty \left\lVert \omega
  \right\rVert _p .
\end{equation*}
To get the lower bound, we begin by using the pointwise inequality
\eqref{eqn:pointwisePullback} for the pullback to write
\begin{equation*}
  \left\lvert \omega \right\rvert = \left\lvert \varphi ^\ast \varphi
    _\ast \omega \right\rvert \leq \alpha _1 \cdots \alpha _k \left(
    \left\lvert \varphi _\ast \omega \right\rvert \circ \varphi
  \right) .
\end{equation*}
Thus, using change of variables and H\"older's inequality once again,
\begin{align*}
  \int _M \left\lvert \omega \right\rvert ^p \mu _M &\leq \int _M
  \left( \alpha _1 \cdots \alpha _k \right) ^p \left( \left\lvert
      \varphi _\ast \omega \right\rvert ^p \circ \varphi \right) \mu
  _M \\
  &= \int _M \varphi ^\ast \left( \left[ \left( \alpha _1 \cdots
        \alpha _k \right) ^p \circ \varphi ^{-1} \right] \left\lvert
      \varphi _\ast \omega \right\rvert ^p \varphi _\ast \mu _M
  \right) \\
  &= \int _M \varphi ^\ast \left( \left[ \left( \alpha _1 \cdots
        \alpha _k \right) ^p \circ \varphi ^{-1} \right] \left\lvert
      \varphi _\ast \omega \right\rvert ^p \bigl[ \left( \alpha _1
      \cdots \alpha _n \right) ^{-1} \circ \varphi ^{-1} \bigr] \mu _N
  \right) \\
  &= \int _M \left( \alpha _1 \cdots \alpha _k \right) ^{ p - 1 }
  \left( \alpha _{ k + 1 } \cdots \alpha _n \right) ^{-1} \varphi
  ^\ast \left( \left\lvert \varphi _\ast \omega \right\rvert ^p \mu _N
  \right) \\
  &\leq \bigl\lVert \mathbf{1} _{ \operatorname{supp} \omega } \left(
    \alpha _1 \cdots \alpha _k \right) ^{ p - 1 } \left( \alpha _{ k +
      1 } \cdots \alpha _n \right) ^{-1} \bigr\rVert _\infty \int _N
  \left\lvert \varphi _\ast \omega \right\rvert ^p \mu _N ,
\end{align*}
so
\begin{equation*}
  \left\lVert \omega \right\rVert _p \leq \bigl\lVert \mathbf{1} _{
    \operatorname{supp} \omega } \left( \alpha _1 \cdots \alpha _k
  \right) ^{ 1/q } \left( \alpha _{ k +  1 } \cdots \alpha _n
  \right) ^{-1/p} \bigr\rVert _\infty \left\lVert \varphi _\ast \omega
  \right\rVert _p .
\end{equation*}
Hence, this implies the lower bound
\begin{equation*}
  \left\lVert \varphi _\ast \omega \right\rVert _p \geq \bigl\lVert
  \mathbf{1} _{ \operatorname{supp} \omega } \left( \alpha _1 \cdots \alpha _k
  \right) ^{ 1/q } \left( \alpha _{ k +  1 } \cdots \alpha _n
  \right) ^{-1/p} \bigr\rVert _\infty  ^{-1}   \left\lVert \omega
  \right\rVert _p .
\end{equation*}
Combining the upper and lower bounds, we have finally established the
two-sided inequality
\begin{multline*}
  \bigl\lVert \mathbf{1} _{ \operatorname{supp} \omega } \left( \alpha
    _1 \cdots \alpha _k \right) ^{ 1/q } \left( \alpha _{ k + 1 }
    \cdots \alpha _n \right) ^{-1/p} \bigr\rVert _\infty ^{-1}
  \left\lVert \omega
  \right\rVert _p  \\
  \leq \left\lVert \varphi _\ast \omega \right\rVert _p \leq
  \bigl\lVert \mathbf{1} _{ \operatorname{supp} \omega } \left( \alpha
    _1 \cdots \alpha _{ n - k } \right) ^{ 1/p } \left( \alpha _{ n -
      k + 1 } \cdots \alpha _n \right) ^{ -1/q } \bigr\rVert _\infty
  \left\lVert \omega \right\rVert _p ,
\end{multline*}
which completes the proof of \autoref{thm:kformPushforward}. \qed

\begin{remark}
  If $\varphi$ is an isometry (at least on $ \operatorname{supp}
  \omega $), then it preserves the $ L ^p $-norm for all $p$.  Indeed,
  isometry implies that the matrix $\Phi$ is orthogonal at every $ x
  \in M $, so the singular values are $ \alpha _1 = \cdots = \alpha _n
  = 1 $.  Therefore, the inequality becomes $ \left\lVert \omega
  \right\rVert _p \leq \left\lVert \varphi _\ast \omega \right\rVert
  _p \leq \left\lVert \omega \right\rVert _p $, and hence $
  \left\lVert \varphi _\ast \omega \right\rVert _p = \left\lVert
    \omega \right\rVert _p $.

  Note that for $ k = 0 $, this simply reduces to
  \autoref{thm:scalarPushforward}, so the $ L ^p $-norm is preserved
  whenever $\varphi$ is volume-preserving (since $ \alpha _1 \cdots
  \alpha _n = 1 $) or when $ p = \infty $, $ q = 1 $.  Similarly, for
  $ k = n $, this reduces to \autoref{thm:densityPushforward}, so the
  $ L ^p $-norm is preserved whenever $\varphi$ is volume-preserving
  or when $ p = 1 $, $ q = \infty $.  More generally, though, for $ 0
  < k < n $, volume preservation is {\em not} sufficient: it merely
  implies that the product of all $n$ singular values equals $1$, but
  it does not imply that this also holds for products of the $k$
  largest or smallest singular values.
\end{remark}

\appendix

\section{Alternating tensors and pointwise inequalities}
\label{app:multilinear}

The pointwise inequalities \eqref{eqn:pointwisePullback} and
\eqref{eqn:pointwisePushforward} are a consequence of some facts from
multilinear algebra regarding alternating tensors.  Consider the case
where $ \Phi $ is a linear isomorphism on $\mathbb{R}^n$ with singular
values $ \alpha _1 \geq \cdots \geq \alpha _n > 0 $.  Associated to
$\Phi$, there is a so-called {\em compound} (or {\em induced
  operator}) on the $k$th exterior power, defined by
\begin{equation*}
  C _k (\Phi) \colon \bigwedge ^k \mathbb{R}^n  \rightarrow
  \bigwedge ^k \mathbb{R}^n  , \quad X _1 \wedge \cdots \wedge X _k
  \mapsto \Phi X _1 \wedge \cdots \wedge \Phi X _k .
\end{equation*}
Due to the total antisymmetry of the exterior product, it can be shown
that $ \left\lvert C _k (\Phi) \right\rvert = \alpha _1 \cdots \alpha
_k $ (see, for example, \citet{AnMa1976,MaAn1977,Li1986}).  This
result is derived by considering the representation of $ C _k (\Phi) $
as an $ \binom{n}{k} \times \binom{n}{k} $ matrix, whose entries are
the corresponding $ k \times k $ minors of $\Phi$, and showing that
its singular values are precisely the $k$-fold products of singular
values of $ \Phi $.  Hence the largest singular value of $ C _k (\Phi)
$ is $ \alpha _1 \cdots \alpha _k $, the product of the $k$ largest
singular values of $\Phi$.  (This basic fact about singular values of
the compound matrix was used at least as early as \citet{Weyl1949},
where it was instrumental in proving the famous Weyl inequalities
relating eigenvalues to singular values.)
It follows that, if $A$ is any alternating $k$-linear form on
$\mathbb{R}^n$, we have
\begin{equation*}
  \bigl( A \circ C _k (\Phi) \bigr) \left( X _1, \ldots , X _k
  \right) = A \left( \Phi X _1 , \ldots, \Phi X _k \right) ,
\end{equation*}
which satisfies the spectral norm inequality
\begin{equation*}
  \left\lvert A \circ C _k (\Phi) \right\rvert \leq \left\lvert A
  \right\rvert \left\lvert C _k (\Phi) \right\rvert = \alpha _1 \cdots
  \alpha _k \left\lvert A \right\rvert .
\end{equation*}

Now, as before, suppose that $ \Phi $ is the matrix representation of
$ T _x \varphi \colon T _x M \rightarrow T _{ \varphi (x) } N $ at
some point $ x \in M $, relative to a positively-oriented, $ g _M
$-orthonormal basis at $x$ and a positively-oriented, $ g _N
$-orthonormal basis at $ \varphi (x) $.  Then the compound matrix $ C
_k (\Phi) $ induces a linear map $ \bigwedge ^k T _x M \rightarrow
\bigwedge ^k T _{ \varphi (x) } N $, and has $ \left\lvert C _k (\Phi)
\right\rvert = \alpha _1 (x) \cdots \alpha _k (x) $.  If $A$ is the
corresponding tensor representation of $ \eta \in \Omega ^k (N) $ at
$\varphi (x) \in N $, then $ \varphi ^\ast \eta $ is given by $ A
\circ C _k (\Phi) $.  Thus,
\begin{equation*}
  \left\lvert \varphi ^\ast \eta \right\rvert = \left\lvert A \circ C
    _k (\Phi) \right\rvert \leq \alpha _1 (x) \cdots \alpha _k (x)
  \left\lvert A \right\rvert =  \alpha _1 (x) \cdots \alpha _k (x)
  \left\lvert \eta \right\rvert \left( \varphi (x) \right)  ,
\end{equation*}
as stated in \eqref{eqn:pointwisePullback}.  Likewise, if $B$
represents $ \omega \in \Omega ^k (M) $ at $ x \in M $, then
\begin{multline*}
  \left\lvert \varphi _\ast \omega \right\rvert = \left\lvert B \circ
    C _k \left( \Phi ^{-1} \right) \right\rvert
  \\
  \leq \beta _1 \left( \varphi (x) \right) \cdots \beta _k \left(
    \varphi (x) \right) \left\lvert B \right\rvert = \beta _1 \left(
    \varphi (x) \right) \cdots \beta _k \left( \varphi (x) \right)
  \left\lvert \omega \right\rvert (x),
\end{multline*}
as in \eqref{eqn:pointwisePushforward}.  (As with the singular values
themselves, these statements are independent of the particular choice
of basis at $ x \in M $ or $ \varphi (x) \in N $.)

Notice that, for $ k = 0 $, there is no contribution from the singular
values, and we just obtain the equalities $ \left\lvert v \circ
  \varphi \right\rvert = \left\lvert v \right\rvert \circ \varphi $
and $ \left\lvert u \circ \varphi ^{-1} \right\rvert = \left\lvert u
\right\rvert \circ \varphi ^{-1} $.  This is why, in the scalar case,
it was possible to simply apply the Jacobian determinant formula to
get
\begin{equation*}
  \varphi ^\ast \left( \left\lvert v \right\rvert ^p \mu _N \right) =
  \left( \left\lvert v \right\rvert ^p \circ \varphi \right) \varphi
  ^\ast \mu _N  = \left\lvert v \circ \varphi \right\rvert ^p \left[ J
    \left( \mu _M , \mu _N \right) \varphi  \right] \mu _M ,
\end{equation*}
while slightly more care was necessary for the density and $k$-form
cases.

\section{Pullback inequalities}
\label{app:pullback}

This appendix contains variants of the main results, stated for
pullbacks rather than pushforwards.  As mentioned in the introduction,
these corollaries follow trivially from the pushforward results given
elsewhere in the paper, simply by replacing $ \varphi $ by $ \varphi
^{-1} $.

\begin{corollary}[smooth functions]
  Let $M$ and $N$ be oriented, $n$-dimensional manifolds with volume
  forms $ \mu _M $ and $ \mu _N $, respectively, and let $ \varphi
  \colon M \rightarrow N $ be an orientation-preserving
  diffeomorphism, whose inverse $ \varphi ^{-1} $ has the Jacobian
  determinant $ J \left( \mu _N , \mu _M \right) \left( \varphi ^{-1}
  \right) \in C ^\infty (N) $.  Then, for any function $ v \in C
  ^\infty (N) $ with compact support $ \operatorname{supp} v $,
  \begin{multline*}
    \bigl\lVert \mathbf{1} _{\operatorname{supp} v } \left[ J \left(
        \mu _N , \mu _M \right) \left( \varphi ^{-1} \right) \right]
    ^{-1/p} \bigr\rVert_\infty ^{-1} \left\lVert v \right\rVert _p \\
    \leq \left\lVert v \circ \varphi \right\rVert _p \leq \bigl\lVert
    \mathbf{1} _{ \operatorname{supp} v } \left[ J \left( \mu _N , \mu
        _M \right) \left( \varphi ^{-1} \right) \right] ^{ 1/p}
    \bigr\rVert _\infty \left\lVert v \right\rVert _p ,
  \end{multline*}
  for all $ p \in \left[ 1, \infty \right] $.
\end{corollary}

\begin{corollary}[$ L ^p $ functions]
  Let $M$ and $N$ be oriented, $n$-dimensional manifolds with volume
  forms $ \mu _M $ and $ \mu _N $, respectively, and let $ \varphi
  \colon M \rightarrow N $ be an orientation-preserving
  diffeomorphism, whose inverse $ \varphi ^{-1} $ has the Jacobian
  determinant $ J \left( \mu _N , \mu _M \right) \left( \varphi ^{-1}
  \right) \in C ^\infty (N) $.  If this inverse Jacobian determinant
  is bounded uniformly on $N$, then for any $ v \in L ^p (N) $,
  \begin{multline*}
    \bigl\lVert \left[ J \left( \mu _N , \mu _M \right) \left( \varphi
        ^{-1} \right) \right] ^{-1/p} \bigr\rVert_\infty ^{-1}
    \left\lVert v \right\rVert _p
    \\
    \leq \left\lVert v \circ \varphi \right\rVert _p \leq \bigl\lVert
    \left[ J \left( \mu _N , \mu _M \right) \left( \varphi ^{-1}
      \right) \right] ^{ 1/p} \bigr\rVert _\infty \left\lVert v
    \right\rVert _p ,
  \end{multline*}
  for all $ p \in \left[ 1, \infty \right] $.
\end{corollary}

\begin{remark}
  By the Jacobian determinant inverse identity
  \eqref{eqn:jacobianInverse}, it follows that $ J \left( \mu _N , \mu
    _M \right) \left( \varphi ^{-1} \right) $ is bounded uniformly on
  $N$ if and only if the reciprocal $ \left[ J \left( \mu _M , \mu _N
    \right) \varphi \right] ^{-1} $ is bounded uniformly on $M$.
  Hence, this inequality can also be written as
  \begin{equation*}
    \bigl\lVert \left[ J \left( \mu _M , \mu _N \right) \varphi
    \right] ^{ 1/p } \bigr\rVert _\infty ^{-1} \left\lVert v
    \right\rVert _p \leq \left\lVert v \circ \varphi \right\rVert _p
    \leq \bigl\lVert \left[ J \left( \mu _M , \mu _N \right) \varphi
    \right] ^{ -1/p } \bigr\rVert _\infty \left\lVert v \right\rVert
    _p .
  \end{equation*}
\end{remark}

\begin{corollary}[smooth densities]
  Let $M$ and $N$ be oriented, $n$-dimensional manifolds with volume
  forms $ \mu _M $ and $ \mu _N $, respectively, and let $ \varphi
  \colon M \rightarrow N $ be an orientation-preserving
  diffeomorphism, whose inverse $ \varphi ^{-1} $ has the Jacobian
  determinant $ J \left( \mu _N , \mu _M \right) \left( \varphi ^{-1}
  \right) \in C ^\infty (N) $. Then, for any smooth density $ v \mu _N
  \in \Omega ^n (N) $, where $ v \in C^\infty (N) $ has compact
  support $ \operatorname{supp} v = \operatorname{supp} v \mu _N $,
  \begin{multline*}
    \bigl\lVert \mathbf{1} _{ \operatorname{supp} v \mu _N }\left[ J
      \left( \mu _N , \mu _M \right) \left( \varphi ^{-1} \right)
    \right] ^{ 1/q } \bigr\rVert _\infty ^{-1}
    \left\lVert v  \mu _N  \right\rVert _p  \\
    \leq \left\lVert \varphi ^\ast \left( v \mu _N \right)
    \right\rVert _p \leq \bigl\lVert \mathbf{1} _{ \operatorname{supp}
      v \mu _N } \left[ J \left( \mu _N , \mu _M \right) \left(
        \varphi ^{-1} \right) \right] ^{ -1/q } \bigr\rVert _\infty
    \left\lVert v \mu _N \right\rVert _p ,
  \end{multline*}
  for all $ p , q \in \left[ 1, \infty \right] $ such that $ 1/p + 1/q
  = 1 $.
\end{corollary}

\begin{corollary}[$ L ^p $ densities]
  Let $M$ and $N$ be oriented, $n$-dimensional manifolds with volume
  forms $ \mu _M $ and $ \mu _N $, respectively, and let $ \varphi
  \colon M \rightarrow N $ be an orientation-preserving
  diffeomorphism, whose inverse $ \varphi ^{-1} $ has the Jacobian
  determinant $ J \left( \mu _N , \mu _M \right) \left( \varphi ^{-1}
  \right) \in C ^\infty (N) $.  If the reciprocal $ \left[ J \left(
      \mu _N , \mu _M \right) \left( \varphi ^{-1} \right) \right]
  ^{-1} \mkern-6mu $ is bounded uniformly on $N$, then for any density
  $ v \mu _N \in L ^p \Omega ^n (N) $, where $ v \in L ^p (N) $,
  \begin{multline*}
    \bigl\lVert \left[ J \left( \mu _N , \mu _M \right) \left( \varphi
        ^{-1} \right) \right] ^{ 1/q } \bigr\rVert _\infty ^{-1}
    \left\lVert v \mu _N \right\rVert _p  \\
    \leq \left\lVert \varphi ^\ast \left( v \mu _N \right)
    \right\rVert _p \leq \bigl\lVert \left[ J \left( \mu _N , \mu _M
      \right) \left( \varphi ^{-1} \right) \right] ^{ -1/q }
    \bigr\rVert _\infty \left\lVert v \mu _N \right\rVert _p ,
  \end{multline*}
  for all $ p , q \in \left[ 1, \infty \right] $ such that $ 1/p + 1/q
  = 1 $.
\end{corollary}

\begin{remark}
  Again, by the Jacobian determinant inverse identity
  \eqref{eqn:jacobianInverse}, the reciprocal $ \left[ J \left( \mu _N
      , \mu _M \right) \left( \varphi ^{-1} \right) \right] ^{-1} $ is
  bounded uniformly on $N$ if and only if $ J \left( \mu _M , \mu _N
  \right) \varphi $ is bounded uniformly on $M$.  Hence, this
  inequality can also be written as
  \begin{multline*}
    \bigl\lVert \left[ J \left( \mu _M , \mu _N \right) \varphi
    \right] ^{ -1/q } \bigr\rVert _\infty ^{-1}
    \left\lVert v \mu _N \right\rVert _p  \\
    \leq \left\lVert \varphi ^\ast \left( v \mu _N \right)
    \right\rVert _p \leq \bigl\lVert \left[ J \left( \mu _M , \mu _N
      \right) \varphi \right] ^{ 1/q } \bigr\rVert _\infty \left\lVert
      v \mu _N \right\rVert _p .
  \end{multline*}
\end{remark}

\begin{corollary}[smooth $k$-forms]
  Let $ \left( M, g _M \right) $ and $ \left( N, g _N \right) $ be
  oriented, $n$-dimensional Riemannian manifolds, and let $ \varphi
  \colon M \rightarrow N $ be an orientation-preserving
  diffeomorphism, whose inverse has singular values $ \beta _1 (y) 
  \geq \cdots \geq \beta _n (y) > 0 $ at each $ y \in N $.  Then, for
  any smooth $k$-form $ \eta \in \Omega ^k (N) $, $ k = 0, \ldots, n
  $, with compact support $ \operatorname{supp} \eta $,
  \begin{multline*}
    \bigl\lVert \mathbf{1} _{\operatorname{supp} \eta } \left( \beta
      _1 \cdots \beta _k \right) ^{ 1/q } \left( \beta _{ k + 1 }
      \cdots \beta _n \right) ^{-1/p} \bigr\rVert _\infty ^{-1}
    \left\lVert \eta
    \right\rVert _p  \\
    \leq \left\lVert \varphi ^\ast \eta \right\rVert _p \leq
    \bigl\lVert \mathbf{1} _{ \operatorname{supp} \eta } \left( \beta
      _1 \cdots \beta _{ n - k } \right) ^{1/p} \left( \beta _{ n - k
        + 1 } \cdots \beta _n \right) ^{ -1/q }\bigr\rVert _\infty
    \left\lVert \eta \right\rVert _p ,
  \end{multline*}
  for all $ p , q \in \left[ 1, \infty \right] $ such that $ 1/p + 1/q
  = 1 $.
\end{corollary}

\begin{corollary}[$ L ^p $ $k$-forms]
  Let $ \left( M, g _M \right) $ and $ \left( N, g _N \right) $ be
  oriented, $n$-dimen\-sional Riemannian manifolds, and let $ \varphi
  \colon M \rightarrow N $ be an orientation-preserving
  diffeomorphism, whose inverse has singular values $ \beta _1 (y)
  \geq \cdots \geq \beta _n (y) > 0 $ at each $ y \in N $.  Given $ p
  , q \in \left[ 1, \infty \right] $ such that $ 1/p + 1/q = 1 $, and
  some $ k = 0 , \ldots, n $, suppose that the product $ \left( \beta
    _1 \cdots \beta _{ n - k } \right) ^{1/p} \left( \beta _{ n - k +
      1 } \cdots \beta _n \right) ^{ -1/q } $ is bounded uniformly on
  $N$.  Then, for any $ \eta \in L ^p \Omega ^k (N) $,
  \begin{multline*}
    \bigl\lVert \left( \beta _1 \cdots \beta _k \right) ^{ 1/q }
    \left( \beta _{ k + 1 } \cdots \beta _n \right) ^{-1/p}
    \bigr\rVert _\infty ^{-1} \left\lVert \eta
    \right\rVert _p  \\
    \leq \left\lVert \varphi ^\ast \eta \right\rVert _p \leq
    \bigl\lVert \left( \beta _1 \cdots \beta _{ n - k } \right)
    ^{1/p} \left( \beta _{ n - k + 1 } \cdots \beta _n \right) ^{
      -1/q }\bigr\rVert _\infty \left\lVert \eta \right\rVert _p .
  \end{multline*}
\end{corollary}

\begin{remark}
  Using the singular value inverse identity $ \beta _i = \alpha _{ n -
    i + 1 } ^{-1} \circ \varphi ^{-1} $, it follows that $ \left(
    \beta _1 \cdots \beta _{ n - k } \right) ^{1/p} \left( \beta _{ n
      - k + 1 } \cdots \beta _n \right) ^{ -1/q } $ is bounded
  uniformly on $N$ if and only if $ \left( \alpha _1 \cdots \alpha _k
  \right) ^{ -1/q } \left( \alpha _{ k + 1 } \cdots \alpha _n \right)
  ^{ 1/p } $ is bounded uniformly on $M$.  Therefore, this inequality
  can also be written as 
  \begin{multline*}
    \bigl\lVert \left( \alpha _1 \cdots \alpha _{n-k} \right) ^{ -1/p
    } \left( \alpha _{ n- k + 1 } \cdots \alpha _n \right) ^{1/q}
    \bigr\rVert _\infty ^{-1} \left\lVert \eta
    \right\rVert _p  \\
    \leq \left\lVert \varphi ^\ast \eta \right\rVert _p \leq
    \bigl\lVert \left( \alpha _1 \cdots \alpha _k \right) ^{-1/q}
    \left( \alpha _{ k + 1 } \cdots \alpha _n \right) ^{ 1/p
    }\bigr\rVert _\infty \left\lVert \eta \right\rVert _p .
  \end{multline*}
\end{remark}

\begin{acknowledgments}
  Thanks to Michael Holst and Tudor Ratiu for reading several earlier
  drafts of this work, and for providing invaluable feedback.
  Research supported in part by NSF (PFC Award 0822283 and DMS Award
  0715146), as well as by NIH, HHMI, CTBP, and NBCR.
\end{acknowledgments}

\end{document}